\documentstyle[12pt,amscd]{amsart}
\newtheorem{thm}[equation]{Theorem}
\numberwithin{equation}{section}

\newtheorem{lem}[equation]{Lemma}

\newtheorem{prop}[equation]{Proposition}

\newtheorem{fig}[equation]{Figure}

\begin{document}
\raggedbottom
\voffset=-.7truein
\hoffset=0truein
\vsize=8truein
\hsize=6truein
\textheight=8truein
\textwidth=6truein
\baselineskip=18truept
\def\mapright#1{\smash{\mathop{\longrightarrow}\limits^{#1}}}
\def\mapleft#1{\smash{\mathop{\longleftarrow}\limits^{#1}}}
\def\mapup#1{\Big\uparrow\rlap{$\vcenter {\hbox {$#1$}}$}}
\def\mapdown#1{\Big\downarrow\rlap{$\vcenter {\hbox {$\ssize{#1}$}}$}}
\def\mapne#1{\nearrow\rlap{$\vcenter {\hbox {$#1$}}$}}
\def\mapse#1{\searrow\rlap{$\vcenter {\hbox {$\ssize{#1}$}}$}}
\def\mapr#1{\smash{\mathop{\rightarrow}\limits^{#1}}}
\def\lb{[}
\def\ss{\smallskip}
\def\sm{\wedge}
\def\la{\langle}
\def\ra{\rangle}
\def\on{\operatorname}
\def\kbar{{\overline k}}
\def\qed{\quad\rule{8pt}{8pt}\bigskip}
\def\ssize{\scriptstyle}
\def\a{\alpha}
\def\bz{\bold Z}
\def\im{\on{im}}
\def\ext{\on{Ext}}
\def\bspin{\on{BSpin}}
\def\sq{\on{Sq}}
\def\eps{\epsilon}
\def\ar#1{\stackrel {#1}{\rightarrow}}
\def\br{\bold R}
\def\bc{\bold C}
\def\si{\sigma}
\def\Ebar{{\overline E}}
\def\Sum{\sum}
\def\tfrac{\textstyle\frac}
\def\tb{\textstyle\binom}
\def\Si{\Sigma}
\def\w{\wedge}
\def\equ{\begin{equation}}
\def\b{\beta}
\def\G{\Gamma}
\def\g{\gamma}
\def\endeq{\end{equation}}
\def\sn{S^{2n+1}}
\def\zp{\bold Z_p}
\def\P{{\cal P}}
\def\zt{\bold Z_2}
\def\Hom{\on{Hom}}
\def\ker{\on{ker}}
\def\coker{\on{coker}}
\def\da{\downarrow}
\def\io{\iota}
\def\Om{\Omega}
\def\u{{\cal U}}
\def\e{{\cal E}}
\def\exp{\on{exp}}
\def\xbar{{\overline x}}
\def\ebar{{\overline e}}
\def\et{{\widetilde E}}
\def\ni{\noindent}
\def\coef{\on{coef}}
\def\den{\on{den}}
\def\lcm{\on{l.c.m.}}
\def\vi{v_1^{-1}}
\def\ot{\otimes}
\def\psibar{{\overline\psi}}
\def\mhat{{\hat m}}
\def\exc{\on{exc}}
\def\ms{\medskip}
\def\ehat{{\hat e}}
\def\BSpt{\widetilde{BSp}}
\def\tbinom#1#2{{\textstyle\binom #1#2}}
\def\dirlim{\on{dirlim}}
\title[Embeddings and immersions]
{Some new embeddings and nonimmersions of real projective spaces}
\author[Davis]{Donald M. Davis}
\address{Lehigh University\\Bethlehem, PA 18015}
\email{dmd1@@lehigh.edu}
\author[Zelov]{Vitaly Zelov}
\address{Lehigh University\\Bethlehem, PA 18015}
\email{viz2@@lehigh.edu}

\date{}

\maketitle

\section{Statement of results} \label{intro}
In this paper, we obtain the following new results regarding immersions and
embeddings of real projective space $P^m$
in Euclidean space. Let $\a(n)$ denote the
number of 1's in the binary expansion of $n$.
\begin{thm} If $\a(n)=2$, then \label{main1}
\begin{enumerate}
\item $P^{16n+8}$ cannot be immersed in $\br^{32n+3}$, and
\item $P^{16n+10}$ cannot be immersed in $\br^{32n+11}$.
\end{enumerate}
\end{thm}
\begin{thm}\label{main2} If $\a(n)>2$, then $P^{8n+4}$ can be embedded in
$\br^{16n+1}$.
\end{thm}

Theorem \ref{main1}(1) improves on the previously best known result (\cite{AD})
by 1 dimension, while Theorem \ref{main1}(2) improves on the previously best
known nonimmersion and nonembedding results (\cite{D83}) for $P^{16n+10}$ and
$P^{16n+11}$ by 4 dimensions, and is within 1 of best possible for them.
It also implies new nonimmersions for
$P^{16n+12}$, $P^{16n+13}$, and $P^{16n+14}$. Theorem \ref{main2} improves
on the previously best known embedding (\cite{Th}) of $P^{8n+4}$ by 1
dimension. These results can best be appreciated when viewed in a table
of known embedding and immersion results for real projective spaces. Such a
table may be seen on the internet at \cite{Dweb}.

The method of proof is obstruction theory, specifically modified Postnikov
towers (MPTs). The reason that Theorem \ref{main1} had not been noticed before
is that the first author used to think that nonimmersions were extremely
difficult to prove by MPTs because of the possibility of secondary and higher
order indeterminacy. We show here that sometimes there is no secondary or
higher order indeterminacy for a simple reason.  The reason that Theorem
\ref{main2} had not been noticed before is that this was apparently the first
time that Mahowald's inductive approach to constructing embeddings
(\cite{Mah64}) was combined
with the method of evaluation of obstructions initiated in
\cite{DM}.

Theorem \ref{main1}(1) and Theorem \ref{main2} form the bulk of the second
author's thesis (\cite{Z}), written under the direction of the first author,
who subsequently discovered the second part of Theorem \ref{main1}.

\section{Proof of nonimmersions}\label{nonimmsect}
In this section, we prove Theorem \ref{main1}(2). The proof of
Theorem \ref{main1}(1), which was detailed in \cite{Z}, is extremely similar,
and is omitted.

The problem is reduced to obstruction theory by the following result of
Sanderson, which reduces the immersion question to the determination of the
geometric dimension of a multiple of the Hopf bundle $\xi_n$ over $P^n$.
\begin{prop} $($\cite{San}$)$ $P^n$ can be immersed in $\br^{n+k}$ if and only
if the map $P^n\to BO$ which classifies $(n+k+1)\xi_n$ can be factored as
$P^n\to BO(k)\to BO$.
\end{prop}
Thus Theorem \ref{main1}(2) will follow from the following result, whose proof
will occupy the rest of this section.
\begin{thm} \label{nonlift} If $\a(n)=2$, then $(32n+12)\xi:P^{16n+10}\to BO$
factors through $BO(16n+1)$.
\end{thm}

The map of Theorem \ref{nonlift} factors as
\begin{equation}\label{map}P^{16n+10}\to HP^{4n+2}@>(8n+3)H>>BSp\to BO.
\end{equation}
Here $HP^{m}$ denotes quaternionic projective space, and $pH$ a multiple of
the quaternionic
Hopf bundle. We let $\BSpt(n)$ denote the classifying space for quaternionic
bundles of real geometric dimension $n$; it is the pullback of maps of
$BO(n)$ and
$BSp$ to $BO$. We will prove Theorem \ref{nonlift} by showing that the
map $P^{16n+10}\to BSp$ in (\ref{map}) lifts to $\BSpt(16n+1)$. This will be
accomplished using the following MPT. For typographical reasons, we abbreviate
$K(\zt,16n+i)$ as $K_i$. We will, however, name the corresponding $k$-invariant
$k^j_{16n+i}\in H^{16n+i}(E_j)$. (All coefficients are in $\bz/2$.) All of our
MPTs are performed through the range of dimensions relevant for the real
projective space being mapped in, here $16n+10$. This will always be well
within the stable range.

\begin{fig}\label{MPT}
\begin{center}\begin{picture}(470,140)
\put(100,0){\makebox(115,15)[br]{$P^{16n+10}\to HP^{4n+2}@>(8n+3)H>>$}}
\put(220,0){$BSp$}
\put(245,0){$\to K_2\times K_4\times K_8$}
\put(230,27){\vector(0,-1){15}}
\put(224,30){$E_1\to K_3\times
K_4\times K_5\times K_7\times K_8\times K_9\times K_9'$}
\put(200,45){\vector(2,-1){20}}
\put(136,41){\makebox(60,15)[br]{$K_1\times K_3\times K_7$}}
\put(230,57){\vector(0,-1){15}}
\put(224,60){$E_2\to K_4\times K_8\times K_8'\times K_9\times K_{10}\times
K_{10}'$}
\put(200,75){\vector(2,-1){20}}
\put(100,71){\makebox(96,15)[br]{$K_2\times K_3\times K_4\times K_6\times
K_7\times K_8\times K_8'$}}
\put(230,87){\vector(0,-1){15}}
\put(224,90){$E_3$}
\put(230,117){\vector(0,-1){15}}
\put(210,120){$\BSpt(16n+1)$}
\end{picture}
\end{center}
\end{fig}

MPTs were introduced in \cite{Mah64} and \cite{GM}.
Each vertical map in the above diagram is
part of a fiber sequence preceded by the map from the fiber represented by a
diagonal arrow, and followed by the classifying map represented by a horizontal
arrow. The information of the diagram can be obtained from the ASS of the
stunted real projective space $P_{16n+1}$, which is, in the stable range,
the fiber of $\BSpt(16n+1)\to BSp$. This ASS can be found in Table 8.2 on
page 54 of \cite{MahMem}. We re-create it in Figure \ref{htpycht}.

 \begin{fig}\label{htpycht}
\begin{center}
\begin{picture}(470,90)
\def\elt{\circle*{3}}
\def\mp{\multiput}
\put(65,5){$*=16n+$}
\put(100,45){$\pi_*(P_{16n+1})$}
\put(138,5){$1$}
\put(167,5){$3$}
\put(227,5){$7$}
\put(130,20){\line(1,0){135}}
\put(140,20){\line(1,1){30}}
\put(170,20){\line(0,1){30}}
\put(230,20){\line(0,1){45}}
\put(230,20){\line(1,1){30}}
\put(247,33){\line(1,1){15}}
\put(245,50){\line(1,1){15}}
\mp(140,20)(15,15){3}{\elt}
\mp(170,20)(0,15){2}{\elt}
\mp(185,35)(30,0){2}{\elt}
\mp(230,20)(0,15){4}{\elt}
\mp(245,35)(15,15){2}{\elt}
\mp(247,33)(15,15){2}{\elt}
\mp(245,50)(15,15){2}{\elt}
\put(233,50){\elt}
\put(260,80){\elt}
\end{picture}
\end{center}
\end{fig}

 The generalized Eilenberg-MacLane spaces (GEMs)
on the right side of Figure \ref{MPT} give rise to the $k$-invariants for the
MPT; these classes have dimensions 1 greater than those of the corresponding
elements of $\pi_*(P_{16n+1})$. We will prove that the map $P^{16n+10}\to
BSp$ lifts to $E_2$, and every lifting to $E_2$ sends some of the level-2
$k$-invariants nontrivially. Hence the map does not lift to $E_3$ or to
$\BSpt(16n+1)$.

The method of evaluation of obstructions is to use the following result
of \cite{DM}. Here $bo$ is the spectrum for connective $ko$-theory localized at
2.
\begin{thm}$($\cite[1.8]{DM}$)$ Let $B^o(m)$ denote the fiberwise smash product
of $\BSpt(m)$ with $bo$. In the stable range, there is a map of fibrations
which is the natural inclusion on the fibers
$$\begin{CD} P_m@>>> P_m\wedge bo\\
@VVV @VVV\\
\BSpt(m)@>>> B^o(m)\\
@VVV @VVV\\
BSp@>=>> BSp
\end{CD}$$
If $pH_k$ denotes $p$ times the Hopf bundle over $HP^k$, then
$pH_k$ lifts to $B^o(m)$ if and only if $m\ge 2k$ and for all $i\le k$
$$\nu(\tbinom {p}{i})\ge\nu(|ko_{4i-1}(P_m)|).$$\label{bothm}
\end{thm}
\noindent Here, and throughout, $\nu(-)$ denotes the exponent of 2.

We use Theorem \ref{bothm} to prove the following result.
\begin{thm}\label{lev2} In the MPT of Figure \ref{MPT}, there is a lifting of
$(8n+3)H_{4n+2}$ to $E_2$ sending $k_{16n+4}^2$ and $k_{16n+8}^2$ nontrivially.
\end{thm}
\begin{pf} We begin by observing that $\nu\binom {8n+3}{4n+1}=\nu\binom{8n+3}
{4n+2}=\a(n)=2$. Indeed, both are equal to $\a(4n+1)+\a(4n+2)-\a(8n+3)$.

These must be compared with the numbers in the following table.
Results such as these have been tabulated in many papers of the first author,
such as \cite{DM}.
\begin{center}
\begin{tabular}{r|cccc|}
&&$m$&&\\
$\nu(ko_{4i-1}(P_m))$&$16n+1$&$16n+2$&$16n+5$&$16n+6$\\
\hline
$i=4n+1$&$3$&$2$&$0$&$0$\\
$i=4n+2$&$4$&$4$&$3$&$2$\\
\hline
\end{tabular}
\end{center}

The orders tabulated here correspond to ASS charts pictured below, which may be
computed, for example, as in \cite{Sur}.
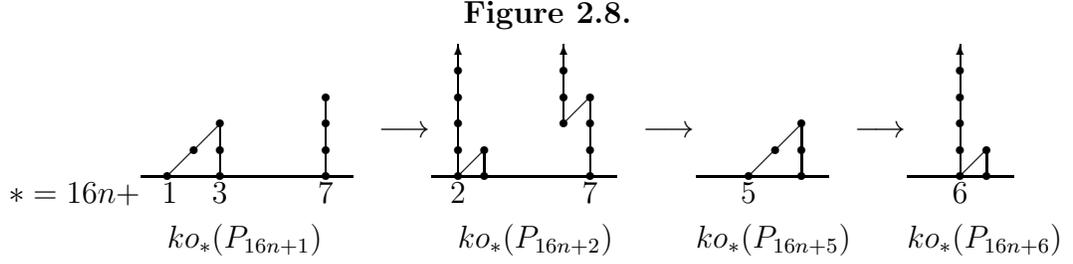
\begin{fig}\label{charts}
\begin{center}
\begin{picture}(480,85)
\def\elt{\circle*{3}}
\def\mp{\multiput}
\put(70,3){$ko_*(P_{16n+1})$}
\put(180,3){$ko_*(P_{16n+2})$}
\put(270,3){$ko_*(P_{16n+5})$}
\put(350,3){$ko_*(P_{16n+6})$}
\put(10,20){$*=16n+$}
\put(68,20){$1$}
\put(87,20){$3$}
\put(127,20){$7$}
\put(177,20){$2$}
\put(227,20){$7$}
\put(287,20){$5$}
\put(367,20){$6$}
\put(60,30){\line(1,0){80}}
\put(170,30){\line(1,0){70}}
\put(270,30){\line(1,0){50}}
\put(350,30){\line(1,0){40}}
\put(70,30){\line(1,1){20}}
\put(90,30){\line(0,1){20}}
\put(130,30){\line(0,1){30}}
\mp(70,30)(10,10){3}{\elt}
\mp(90,30)(0,10){2}{\elt}
\mp(130,30)(0,10){4}{\elt}
\put(150,45){$\longrightarrow$}
\put(250,45){$\longrightarrow$}
\put(330,45){$\longrightarrow$}
\put(180,30){\vector(0,1){50}}
\put(180,30){\line(1,1){10}}
\put(190,30){\line(0,1){10}}
\put(220,50){\vector(0,1){30}}
\put(220,50){\line(1,1){10}}
\put(230,30){\line(0,1){30}}
\mp(180,30)(0,10){5}{\elt}
\mp(190,30)(0,10){2}{\elt}
\mp(220,50)(0,10){3}{\elt}
\mp(230,30)(0,10){4}{\elt}
\put(290,30){\line(1,1){20}}
\put(310,30){\line(0,1){20}}
\mp(290,30)(10,10){3}{\elt}
\mp(310,30)(0,10){2}{\elt}
\put(370,30){\vector(0,1){50}}
\put(370,30){\line(1,1){10}}
\put(380,30){\line(0,1){10}}
\mp(370,30)(0,10){5}{\elt}
\mp(380,30)(0,10){2}{\elt}
\end{picture}
\end{center}
\end{fig}

There are MPTs for the fibrations $B^o(16n+\delta)\to BSp$, and induced maps of
the MPTs, whose spaces we denote by $E_s^o(16n+\delta)$.
By Theorem \ref{bothm}, $(8n+3)H_{4n+1}$ lifts to $B^o(16n+2)$ but not to
$B^o(16n+1)$. This implies that $(8n+3)H_{4n+1}$ lifts to $E_2^o(16n+1)$
sending $k^2_{16n+4}$ nontrivially. (Keep in mind that degrees of
$k$-invariants are 1 greater than those of corresponding elements of
$\pi_*(\text{fiber})$, which in this case is $ko_*(P_{16n+\delta})$.)
Similarly, $(8n+3)H_{4n+2}$ lifts to $B^o(16n+6)$, but not to
$B^o(16n+5)$. All this implies that $(8n+3)H_{4n+2}$ lifts to $E_2^o(16n+1)$,
sending both $k^2_{16n+4}$ and $k^2_{16n+8}$ nontrivially.

The map $\BSpt(16n+1)\to B^o(16n+1)$ induces a map of MPTs. The mapping of
$k$-invariants corresponds to the map of fibers $\pi_*(P_{16n+1})\to
ko_*(P_{16n+1})$, which are pictured in Figures \ref{htpycht} and \ref{charts}.
Various methods of computing Ext homomorphisms can be used to show this
morphism is surjective in the indicated range. If $F$ denotes the fiber of
$E_2\to E_2^o(16n+1)$, then $\pi_*(F)$ corresponds to elements in filtration 0
and 1 in Figure \ref{htpycht} which map trivially to $ko_*(P_{16n+1})$.
Such elements occur only in filtration 1 and in homotopy dimension $16n+4$,
$16n+6$, and $16n+8$. The obstructions to pulling the map $HP^{4n+2}\to
E_2^o(16n+1)$ back to $E_2$ occur in $H^*(HP^{4n+2};\pi_{*-1}(F))$, which is 0
since $\pi_{*-1}(F)=0$ when $*\equiv0$ mod 4
and $*\le16n+8$. Thus $(8n+3)H_{4n+2}$ lifts to
$E_2$. The $k$-invariants $k^2_{16n+4}$ and $k^2_{16n+8}$ in $H^*(E_2)$
are the images of corresponding $k$-invariants in $H^*(E_2^o(16n+1))$, which
have already been shown to map nontrivially to $H^*(HP^{4n+2})$. Thus
$k^2_{16n+4}$ and $k^2_{16n+8}$ map nontrivially, as claimed. We make no claim
about whether or not $k^{2'}_{16n+8}$, which corresponds to the split $\bz/2$
in $\pi_{16n+7}(P_{16n+1})$ in filtration 2, maps nontrivially, since it is not
in the image from $H^*(E_2^o(16n+1))$.
\end{pf}

In order to determine the indeterminacy for lifting $P^{16n+10}$ in this MPT,
we must know the relations which give rise to the $k$-invariants. These are
computed by the method initiated in \cite{GM} and utilized in many subsequent
papers by the first author and also in papers of Lam and/or Randall.
It is a matter of building a minimal resolution using Massey-Peterson algebras.
The relations for the MPT in Figure \ref{MPT} are given in the table below.

\begin{center}
\begin{tabular}{|ll|}
\hline
$w_{16n+2}$&\\
$w_{16n+4}$&\\ $w_{16n+8}$&\\
\hline
$k^1_{16n+3}:$&$\ \sq^2w_{16n+2} $\\
$k^1_{16n+4}:$&$\ \sq^1w_{16n+4}+\sq^2\sq^1w_{16n+2}$\\
$k^1_{16n+5}:$&$\ (\sq^4+w_4)w_{16n+2}$\\
$k^1_{16n+7}:$&$\ (\sq^4+w_4)w_{16n+4}$\\
$k^1_{16n+8}:$&$\ \sq^1w_{16n+8}+\sq^2\sq^3w_{16n+4}$\\
$k^1_{16n+9}:$&$\ \sq^2w_{16n+8}+(\sq^4+w_4)\sq^2w_{16n+4}$\\
$k^{1'}_{16n+9}:$&$\ (\sq^8+w_8)w_{16n+2}+w_4\sq^2w_{16n+4}$\\
\hline
$k^2_{16n+4}:$&$\ \sq^2k^1_{16n+3}+\sq^1k^1_{16n+4}$\\
$k^2_{16n+8}:$&$\ \sq^1k^1_{16n+8}+\sq^2\sq^3k^1_{16n+4}$\\
$k^{2'}_{16n+8}:$&$\
(\sq^4+\sq^3\sq^1+w_4)k^1_{16n+5}+(\sq^6+w_4\sq^2)k^1_{16n+3}$\\
$k^2_{16n+9}:$&$\ \sq^2\sq^1k^1_{16+7}+(\sq^4+w_4)\sq^1k^1_{16n+5}+(\sq^6+
w_4\sq^2)k^1_{16n+4}$\\
&$\ +(\sq^4w_4)\sq^2\sq^1k^1_{16n+3}$\\
$k^2_{16n+10}:$&$\
\sq^2k^1_{16n+9}+\sq^3k^1_{16n+8}+(\sq^4+\sq^3\sq^1+w_4)k^1_{16n+7}$\\
$k^{2'}_{16n+10}:$&$\ \sq^2k^{1'}_{16n+9}+(\sq^5\sq^1+\sq^4\sq^2)k^1_{16n+5}
+w_4\sq^3k^1_{16n+4}$\\
&$\ +(\sq^8+w_8+w_4\sq^4+w_4^2+w_4\sq^3\sq^1)k^1_{16n+3}$\\
\hline
$k^3_{16n+8}:$&$\ \sq^1k^2_{16n+8}+\sq^2\sq^3k^2_{16n+4}$\\
$k^3_{16n+10}:$&$\
\sq^2k^2_{16n+9}+\sq^2\sq^1k^{2'}_{16n+8}$\\
&$\ +(\sq^7+\sq^4\sq^2\sq^1+w_4(\sq^3+\sq^
2\sq^1)k^2_{16n+4}$\\
\hline
$k^4_{16n+10}:$&$\ \sq^2\sq^1k^3_{16n+8}$\\
\hline
\end{tabular}
\end{center}

The lifting $f_2:P^{16n+10}\to E_2$ can be varied through the fiber $F_1$ of $E_
2\to
E_1$, which is the GEM in Figure \ref{MPT} which ends with $K_8'$.
This primary indeterminacy
is computed using the above relations. For the bundle $(32n+12)\xi$, both
$w_4$ and $w_8$ are nonzero. Varying through $K_{16n+2}$  changes
$f_2^*(k^2_{16n+4})$ and $f_2^*(k^{2'}_{16n+8})$ and no other level-2
$k$-invariants. We illustrate why this is
true in the second case.

The relation for $k^{2'}_{16n+8}$ means that the action map
$\mu:F_1\times E_2\to E_2$ sends $k^{2'}_{16n+8}$ to
$$1\ot k^{2'}_{16n+8}+(\sq^4+\sq^3\sq^1)\io_{16n+4}\ot1+\io_{16n+4}\ot w_4
+\sq^6\io_{16n+2}\ot1+\sq^2\io_{16n+2}\ot w_4.$$
The composite
\begin{equation}\label{comp}
P^{16n+10}@>x^{16n+2}\times f_2>> F_1\times E_2@>\mu>> E_2\end{equation}
sends $k^{2'}_{16n+8}$ to
$$f_2^*(k^{2'}_{16n+8})+\sq^6(x^{16n+2})+\sq^2(x^{16n+2})\cdot x^4=
f_2^*(k^{2'}_{16n+8})+0+x^{16n+8}.$$
Thus (\ref{comp}) would be a new lifting to $E_2$ with $f_2^*(k^{2'}_{16n+8})$
changed. A similar computation is required to
determine whether (\ref{comp}) changes each
of the other $k$-invariants, and the result is as claimed in the previous
paragraph.

Similarly, varying through $K_{16n+3}$ changes $f_2^*(k^2_{16n+4})$,
$f_2^*(k^2_{16n+8})$, $f_2^*(k^2_{16n+9})$, and $f_2^*(k^{2'}_{16n+10})$.
These occur because the following terms, respectively, are nonzero:
$\sq^1x^{16n+3}$, $\sq^2\sq^3(x^{16n+3})$, $w_4\sq^2(x^{16n+3})$, and
$w_4\sq^3(x^{16n+3})$.
The only other way in which $f_2^*(k^2_{16n+\delta})$
can be changed is by varying through
$K_{16n+7}$, which changes $f_2^*(k^2_{16n+8})$ and $f_2^*(k^2_{16n+10})$.

Thus varying the lifting through $F_1$ in such a way that both
$f_2^*(k_{16n+4}^2)$ and $f_2^*(k^2_{16n+8})$ are both changed to become 0
will cause either $f_2^*(k^2_{16n+10})$ or $f_2^*(k^{2'}_{16n+10})$ to become
nonzero. Thus if $f_1:P^{16n+10}\to E_1$ represents a lifting obtained by
factoring through $HP^{4n+2}$, then every lifting of $f_1$ to $E_2$ sends some
$k$-invariants nontrivially, and hence $f_1$ does not lift to $\BSpt(16n+1)$.

We must also consider the possibility that $f_1$ could be varied through
the fiber, $F_0$, of the map $E_1\to BSp$ in such a way that the new map $f_1'$
lifted to $f_2':P^{16n+10}\to E_2$ and sent all $k^2$-invariants to 0. This
is the secondary indeterminacy consideration that had led the first author
to not try to prove nonimmersions by ordinary MPTs in his work during the 1970s
and 1980s.  However, in the case at hand, secondary
indeterminacy is not a problem because any nontrivial map $P^{16n+10}\to F_0$
will change the images of some of the $k^1$-invariants. Since $f_1$ was a map
which lifted to $E_2$ and therefore sent all the $k^1$-invariants to 0, any
variation of $f_1$ through $F_0$ will not lift to $E_2$.

This completes the proof that $(32n+12)\xi_{16n+10}$ does not lift to
$\BSpt(16n+1)$ when $\a(n)=2$, once we verify the statement made in the
preceding paragraph about varying maps through $F_0$. Varying through
$K_{16n+1}$ changes $f_1^*(k^1_{16n+4})$, $f_1^*(k^1_{16n+5})$, and $f_1^*(
 k^{1'}_{16n+9})$. Varying through $K_{16n+3}$ changes $f_1^*(k^1_{16n+4})$,
$f_1^*(k_{16n+7}^1)$, $f_1^*(k^1_{16n+8})$, and $f_1^*(k^{1'}_{16n+9})$.
Varying through $K_{16n+7}$ changes $f_1^*(k_{16n+8}^1)$ and
$f_1^*(k_{16n+9}^1)$. Any nontrivial combination of these changes some
$f_1^*(k^1_{16n+\delta})$, as claimed.

\section{Proof of embeddings}\label{embeddings}
In this section, we prove Theorem \ref{main2}. We use the following result of
Mahowald, which deals with topological embeddings.

\begin{thm} $($\cite{Mah64}$)$
Assume that $P^q$ embeds in $\br^p$ with normal bundle $\nu$.
\begin{itemize}
\item If $\nu\otimes\xi_q$ has $n$ linearly independent sections
and $P^{n-1}$ embeds in
$S^{m-1}$, then $P^{n+q}$ embeds in $\br^{p+m}$.
\item $\nu\otimes\xi_q\oplus(q+1)\eps\approx(p+1)\xi_q$.
\end{itemize}\label{tool}\end{thm}

We apply this result to the embedding of $P^{8n+2}$ in $\br^{16n-1}$
when $\a(n)>2$ proved in
\cite{Th}. Using also the embedding of $P^1$ in $S^1$, we obtain Theorem
\ref{main2}
once we prove the following result. Here $\theta=\nu\ot\xi_{8n+2}$.
\begin{thm}\label{work} If $\theta$ is an $(8n-3)$-plane bundle over $P^{8n+2}$
which is stably equivalent to $16n\xi_{8n+2}$, and $\a(n)>2$, then $\theta$
has at least $2$ linearly independent sections.
\end{thm}
Theorem \ref{tool} will then imply that there is a topological embedding of
$P^{8n+4}$ in $\br^{16n+1}$. Such an embedding can be approximated by a
differentiable embedding by a result of Haefliger (\cite{Haef}), since
$2(16n+1)\ge 3(8n+4)$.

The first step toward proving Theorem \ref{work} is to prove the following
lemma.
\begin{lem}\label{stable} If $\a(n)>2$, the map $P^{8n+2}\to BSp$ which
classifies  $16n\xi_{8n+2}$ lifts to $\BSpt(8n-5)$.
\end{lem}
\begin{pf} Similarly to the work of the previous section, we consider the
following diagram.
$$\begin{array}{rccc}
&P_{8n-5}&\mapright{i}&P_{8n-5}\wedge bo\\
&\da&&\da\\
&\BSpt(8n-5)&\mapright{j}&B^o(8n-5)\\
&\da&&\da\\
P^{8n+2}\to HP^{2n}@>4nH>>&BSp&\mapright{=}&BSp
\end{array}$$

The morphism $\pi_*(P_{8n-5})\to ko_*(P_{8n-5})$ of homotopy groups of fibers
is depicted in the following ASS charts. Again, the first is from \cite{MahMem}
while the second is well known, e.g. \cite{Sur}. The morphism is easily seen to
be surjective.
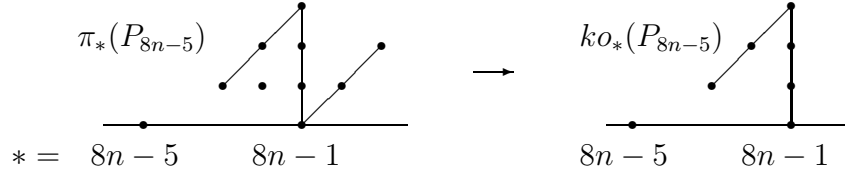
\begin{center}
\begin{picture}(310,70)
\def\mp{\multiput}
\def\elt{\circle*{3}}
\put(-5,5){$*=$}
\put(25,5){$8n-5$}
\put(86,5){$8n-1$}
\put(210,5){$8n-5$}
\put(271,5){$8n-1$}
\put(30,20){\line(1,0){115}}
\put(220,20){\line(1,0){90}}
\put(45,20){\elt}
\mp(75,35)(15,15){3}{\elt}
\mp(105,20)(0,15){3}{\elt}
\mp(120,35)(15,15){2}{\elt}
\put(90,35){\elt}
\put(75,35){\line(1,1){30}}
\put(105,20){\line(0,1){45}}
\put(105,20){\line(1,1){30}}
\put(170,40){\vector(1,0){15}}
\put(230,20){\elt}
\put(260,35){\line(1,1){30}}
\put(290,20){\line(0,1){45}}
\mp(260,35)(15,15){3}{\elt}
\mp(290,20)(0,15){3}{\elt}
\put(20,50){$\pi_*(P_{8n-5})$}
\put(210,50){$ko_*(P_{8n-5})$}
\end{picture}
\end{center}

From the chart, we see that $\nu(|ko_{8n-5}(P_{8n-5})|)=1$, and $\nu(|ko_{8n-1}
(P_{8n-5})|)=4$. One easily calculates $\nu(\binom{4n}{2n-1})>2$ and
$\nu(\binom{4n}{2n})=\a(n)$. If $\a(n)>3$, then, by Theorem \ref{bothm},
$4nH_{2n}$ lifts to a map $\ell:HP^{2n}\to B^o(8n-5)$. Let
$F=\text{fiber}(j)=\text{fiber}(i)$. Since $\pi_*(F)=0$ when $*\equiv3$ mod 4
and $*\le 8n+2$,
$\ell$ pulls back to the desired map $HP^{2n}\to \BSpt(8n-5)$.

When $\a(n)=3$, we can show, similarly to the proof of Theorem \ref{lev2},
that $4nH_{2n}$ lifts to $E_3^o(8n-5)$, the third stage of the MPT for $B^o(8n-
5)$. This uses Theorem \ref{bothm} to see that $4nH_{2n}$ lifts to $B^o(8n-3)$
and $4nH_{2n-1}$ lifts to $B^o(8n-5)$, which imply that the obstructions to
lifting to $E_3^o(8n-5)$ map trivially. The map $HP^{2n}\to E_3^o(8n-5)$ pulls
back to $E_3$, the third stage in the MPT for $\BSpt(8n-5)$, since
$H^*(HP^{2n};\pi_{*-1}(\text{fiber}))=0$.

The fiber of $E_3\to E_2$ contains a $K_{8n-1}$-factor, corresponding to the
third dot up in the ASS chart for $\pi_{8n-1}(P_{8n-5})$.
The map $f_3:RP^{8n+2}\to HP^{2n}\to E_3$ can be varied through this factor
to change $f_3^*(k^3_{8n})$. This is implied by a $\sq^1k^2_{8n}$ which
appears in the relation for $k^3_{8n}$ in the appropriate MPT.
Since $k^3_{8n}$ is the only level-3 $k$-invariant, we deduce that
there is a choice of the lifting $f_3$ which lifts to $E_4=\BSpt(8n-5)$.
\end{pf}

The bundle $\theta$ with which we ultimately must deal might not be a
symplectic bundle; however, it certainly is a Spin bundle. We reinterpret Lemma
\ref{stable} to say that the map $16n\xi:P^{8n+2}\to \bspin$ lifts to
$\bspin(8n-5)$.
Next we consider the following diagram.
\begin{eqnarray}
V_{8n-3,2}&\to&P_{8n-5}\nonumber\\
\da\quad&&\quad\da\label{lastdiag}\\
\bspin(8n-5)&\mapright{=}&\bspin(8n-5)\nonumber\\
\da\quad&&\quad\da\nonumber\\
P^{8n+2}\mapright{\theta}\bspin(8n-3)&\to&\bspin\nonumber
\end{eqnarray}

The first fiber, $V_{8n-3,2}$, is a Stiefel manifold, and in the stable range
is homotopy equivalent to the stunted real projective space $P_{8n-5}^{8n-4}$.
The induced morphism of homotopy groups of fibers can be easily determined
using the tables of \cite{MahMem} to be as below, where the big dots map
across. It is useful to use here
that these two charts fit into an exact sequence in which the third chart
 is $\pi_*(P_{8n-3})$.
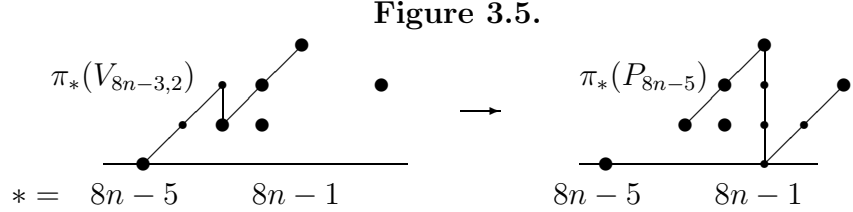
\begin{fig}\label{bigelts}
\begin{center}
\begin{picture}(460,70)
\def\mp{\multiput}
\def\elt{\circle*{3}}
\def\bigelt{\circle*{5}}
\put(45,5){$*=$}
\put(75,5){$8n-5$}
\put(136,5){$8n-1$}
\put(250,5){$8n-5$}
\put(311,5){$8n-1$}
\put(95,20){\bigelt}
\put(95,20){\line(1,1){30}}
\mp(110,35)(15,15){2}{\elt}
\put(125,35){\line(0,1){15}}
\put(125,35){\line(1,1){30}}
\mp(125,35)(15,15){3}{\bigelt}
\put(140,35){\bigelt}
\put(185,50){\bigelt}
\put(80,20){\line(1,0){115}}
\put(260,20){\line(1,0){90}}
\put(270,20){\bigelt}
\mp(300,35)(15,15){3}{\bigelt}
\mp(330,20)(0,15){3}{\elt}
\put(345,35){\elt}
\put(360,50){\bigelt}
\put(315,35){\bigelt}
\put(300,35){\line(1,1){30}}
\put(330,20){\line(0,1){45}}
\put(330,20){\line(1,1){30}}
\put(215,40){\vector(1,0){15}}
\put(60,50){$\pi_*(V_{8n-3,2})$}
\put(260,50){$\pi_*(P_{8n-5})$}
\end{picture}
\end{center}
\end{fig}

The MPT for the fibration $\bspin(8n-5)\to\bspin(8n-3)$ corresponds to the
above chart of $\pi_*(V_{8n-3,2})$. We denote its spaces by $A_i$.

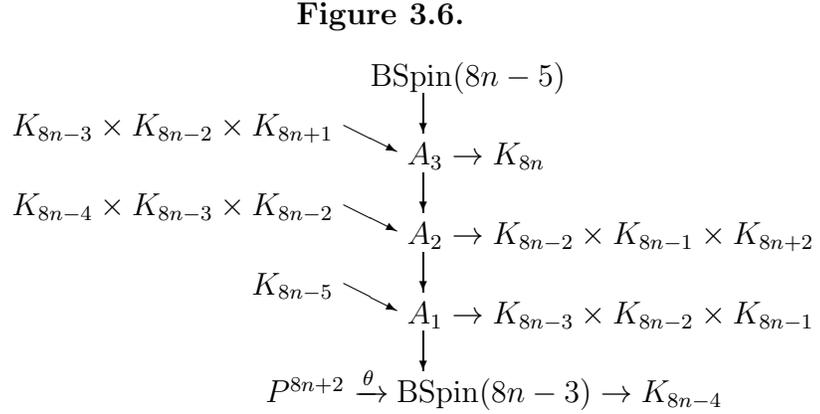
\begin{fig}\label{MPT2}
\begin{center}\begin{picture}(470,140)
\put(100,0){\makebox(115,15)[br]{$P^{8n+2}@>\theta>>$}}
\put(220,0){$\bspin(8n-3)\to K_{8n-4}$}
\put(230,27){\vector(0,-1){15}}
\put(224,30){$A_1\to K_{8n-3}\times K_{8n-2}\times K_{8n-1}$}
\put(200,45){\vector(2,-1){20}}
\put(136,41){\makebox(60,15)[br]{$K_{8n-5}$}}
\put(230,57){\vector(0,-1){15}}
\put(224,60){$A_2\to K_{8n-2}\times K_{8n-1}\times K_{8n+2}$}
\put(200,75){\vector(2,-1){20}}
\put(200,105){\vector(2,-1){20}}
\put(100,71){\makebox(96,15)[br]{$K_{8n-4}\times K_{8n-3}\times K_{8n-2}$}}
\put(230,87){\vector(0,-1){15}}
\put(224,90){$A_3\to K_{8n}$}
\put(230,117){\vector(0,-1){15}}
\put(210,120){$\bspin(8n-5)$}
\put(100,101){\makebox(96,15)[br]{$K_{8n-3}\times K_{8n-2}\times K_{8n+1}$}}
\end{picture}
\end{center}
\end{fig}

We need to know the relations that give rise to these $k$-invariants.
We list them without listing the $w_4$ and $w_8$ which appear in some of these
relations, because $w_4=0$ and $w_8=0$ for the bundle being considered here.

\begin{center}
\begin{tabular}{|l|}
\hline
$w_{8n-4}$\\
\hline
$k^1_{8n-3}:\sq^2w_{8n-4}$\\
$k^1_{8n-2}:\sq^2\sq^1w_{8n-4}$\\
$k^1_{8n-1}:\sq^4w_{8n-4}$\\
\hline
$k^2_{8n-2}:\sq^1k^1_{8n-2}+\sq^2k^1_{8n-3}$\\
$k^2_{8n-1}:\sq^2k^1_{8n-2}$\\
$k^2_{8n+2}:(\sq^4+\sq^3\sq^1)k_{8n-1}^1+\sq^6k_{8n-3}^1$\\
\hline
$k^3_{8n}:\sq^3k^2_{8n-2}+\sq^2k_{8n-1}^2$\\
\hline
\end{tabular}
\end{center}

Primary indeterminacy allows the following classes to be changed (if liftings
to the indicated level exist):
\begin{itemize}
\item $f_1^*(k^1_{8n-3})$ (through $K_{8n-5}$);
\item $f_2^*(k^2_{8n-2})$ (through $K_{8n-3}$);
\item $f_2^*(k^2_{8n+2})$ (through $K_{8n+2}$);
\item $f_3^*(k^3_{8n})$ (through $K_{8n}$).
\end{itemize}
We can also use relations in the MPT to deduce that (again, if the liftings to
the indicated level exist) $f_1^*(k^1_{8n-2})=0$ and $f_2^*(k^2_{8n-
1})=0$. To see the first of these, we use the relation that gives rise to
$k^2_{8n-1}$. This says that $\sq^2k^1_{8n-2}=0\in H^*A_1$. Thus $\sq^2(f_1^*(
k^1_{8n-2}))=0\in H^*(P^{8n+2})$, and so we cannot have $f_1^*(k^1_{8n-2})=
x^{8n-2}$. The equation $f_2^*(k^2_{8n-1})=0\in H^*A_2$ follows similarly
from the relation for $k^3_{8n}$.

Because $w_{8n-4}(\theta)=w_{8n-4}(16n\xi)=\binom{16n}{8n-4}x^{8n-4}=0\in
H^*(P^{8n+2})$, the map $\theta$ lifts to $f_1:P^{8n+2}\to A_1$.
Once we show that $f_1^*(k^1_{8n-1})=0$, we will be able to deduce our desired
lifting to $\bspin(8n-5)$, for all other $k$-invariants are either in the
indeterminacy or map to 0 by relations, as described in the preceding
paragraph.

The map of fibrations in (\ref{lastdiag}) induces a map of MPTs. The element
in filtration 1 in $\pi_{8n-2}(\text{fibers})$ maps across, as indicated in
Figure
\ref{bigelts}. This implies that the corresponding level-1 $k$-invariant
maps across; i.e., if $g_1:A_1\to E_1$ is the induced map of spaces in the
MPTs, then $g_1^*(\tilde k^1_{8n-1})=k^1_{8n-1}$, where we use the tilde to
indicate the $k$-invariants in the MPT for $\bspin(8n-5)\to\bspin$.
In the paragraph after the proof of Lemma \ref{stable}, we have observed that
the map $\theta$ into $\bspin$ lifts to $\bspin(8n-5)$. This implies that
the composite $P^{8n-2}\mapright{f_1} A_1\mapright{g_1}E_1$ differs from a map
$h_1$ that lifts to $\bspin(8n-5)$ by a map $\delta_1$
which factors through the fiber $K_{8n-5}\times K_{8n-1}$
of $E_1\to\bspin$. The relation for
$\tilde k^1_{8n-1}$ is $(\sq^4+w_4)w_{8n-4}$. Since $\sq^4(x^{8n-5})=0$ and
$w_4(\theta)=0$,
any such
map $\delta_1$ must send $\tilde k^1_{8n-1}$ to 0. Since $h_1$ lifts, it
sends all $k$-invariants to 0. Thus $0
=h_1^*(\tilde k^1_{8n-1})+\delta_1^*(\tilde k^1_{8n-1})=
f_1^*g_1^*(\tilde k^1_{8n-1})=
f_1^*(k^1_{8n-1})$, as desired.

\def\line{\rule{.6in}{.6pt}}

\end{document}